\documentclass{article}
\usepackage{latexsym,amssymb,amsmath}
\usepackage{graphics,amsmath,amssymb}
\usepackage{amsthm}
\usepackage{amsfonts}
\usepackage{latexsym}
\usepackage{epsf}
\newtheorem{theorem}{Theorem}
\newtheorem{corollary}{Corollary}
\newtheorem{note}{Note}

\title{On  Non-central Stirling Numbers of the First Kind}
\author{Milan Janji\'c}
\date{}
\begin{document}
\maketitle
\begin{abstract}
It is shown in this note that non-central Stirling numbers $s(n,k,\alpha)$ of the first kind  naturally appear in the expansion of derivatives of the function $x^{-\alpha}\ln^\beta x,$ where $\alpha$ and $\beta$ are arbitrary real numbers. We first obtain a recurrence relation for these numbers, and then, using Leibnitz rule we obtain an explicit formula for them.  We also obtain a formula for $s(n,1,\alpha)$ and then derive several combinatorial identities related to these numbers.
\end{abstract}

\section{ Introduction}
We are dealing here with a special kind of numbers introduced by D. S. Mitrinovi\'c in his note [4].
In the paper [5] tables are given for the numbers which we called non-central Stirling numbers of the first kind.  Following [3], they will be denoted by $s(n,k,\alpha).$  Several other names are in use for these numbers. One of them is $r$-Stirling numbers, as in [1]. The definition in this paper is restricted to the case when $\alpha$ is an nonnegative integer, and $\alpha\leq n.$ L.  Carlitz [2] used the name: weighted Stirling numbers. In the well known encyclopedia [6] they are called the generalized Stirling numbers.
Here we use the name and the notation from the book [3]. For instance,
 $(\alpha)_n=\alpha(\alpha-1)\cdots(\alpha-n+1)$ are
falling factorials, and $s(n,k)$ are Stirling numbers of the first kind.

\section{Derivatives of $x^{-\alpha}\ln^\beta x.$}
We shall investigate derivatives of the function   \[f(x)=x^{-\alpha}\ln^\beta x,\;(\alpha,\beta\in \mathbb R),\]
 obtaining them in two different way.

\begin{theorem} Let $\alpha$ be real, and $n$ nonnegative integer. Then
\begin{equation}\label{e1}f^{(n)}(x)=x^{-\alpha-n}\sum_{i=0}^ns(n,i,\alpha)(\beta)_i\ln^{\beta-i} x.\end{equation}
where $s(n,i,\alpha),\;(0\leq i\leq n)$ are polynomials of $\alpha$ with integer coefficients.
\end{theorem}
 Proof.
The assertion is true for $n=0$ if we take
$s(0,0,\alpha)=1.$

Taking
\[s(1,0,\alpha)=-\alpha,\;s(1,1,\alpha)=1,\]
we see that the assertion is true  for $n=1.$

Suppose that the assertion is true for $n\geq 1.$

Taking derivative in (\ref{e1}) we obtain
\[f^{(n+1)}(x)=(-\alpha-n)x^{-\alpha-n-1}\sum_{i=0}^ns(n,i,\alpha)(\beta)_i\ln^{\beta-i}x+\]\[ +x^{-\alpha-n-1}\sum_{i=0}^{n-1}s(n,i,\alpha)(\beta)_{i+1}\ln^{\beta-i-1}x=\]
\[=x^{-\alpha-n-1}\left[(-\alpha-n)\sum_{i=0}^ns(n,i,\alpha)(\beta)_i\ln^{\beta-i}x+\sum_{i=0}^{n} s(n,i,\alpha)(\beta)_{i+1}\ln^{\beta-i-1}x.\right]\]
Replacing $i+1$ by $i$ in the second sum we obtain
\[f^{(n+1)}(x)=x^{-\alpha-n-1}\left[\sum_{i=0}^n(-\alpha-n)s(n,i,\alpha)(\beta)_i+\sum_{i=1}^{n+1}s(n,i-1,\alpha)(\beta)_{i}\ln^{\beta-i}x\right]=\]
\[=x^{\alpha-n-1}(\alpha-n)s(n,0,\alpha)+\]\[+x^{\alpha-n-1}\sum_{i=1}^n\left[(-\alpha-n)s(n,i,\alpha)+s(n,i-1,\alpha)\ln^{\beta-i}x\right](\beta)_i+\]
\[+s(n,n,\alpha)(\beta)_{n+1}\ln^{\beta-n-1}x.\]
It follows that the  assertion is true if we take
\[ s(n+1,0,\alpha)=(-\alpha-n)s(n,0,\alpha),\]
\begin{equation}\label{e2}s(n+1,i,\alpha)=(-\alpha-n)s(n,i,\alpha)+s(n,i-1,\alpha),\;(i=1,\ldots,n),\end{equation}
\[s(n+1,n+1,\alpha)=s(n,n,\alpha).\]

The preceding equation are well-known recurrence relations for non-central Stirling numbers of the first kind [3, pp.316].
\begin{note}
It is obvious that $s(n,i,0)=s(n,i))$ are Stirling numbers of the first kind.
\end{note}
Since
$p(0,0,\alpha)=1,$ for from $n=1,2,\ldots$ from the first equation in (\ref{e2}) we obtain:
\[p(n,0,\alpha)=(-\alpha)_n,\]
and, since $p(1,1,\alpha)=1,$  from the last equation in (\ref{e2}) follows:
\[p(n,n,\alpha)=1,\;(n=0,1,2,\ldots).\]

By the use of Leibnitz formula we shall obtain an explicit expression for $s(n,k,\alpha).$
The following equation holds:
\begin{equation}\label{e3}f^{(n)}(x)=\sum_{k=0}^n{n\choose k}\big(x^\alpha\big)^{(k)}(\ln^\beta x)^{(n-k)}.\end{equation}

First we have
\[\big(x^{-\alpha}\big)^{(k)}=(-\alpha)_kx^{\alpha-k}.\]

Using induction it is easy to prove that:
\[[f(\ln x)]^{(n)}=x^{-n}\sum_{k=1}^ns(n,k)f^{(k)}(t),\;(t=\ln x).\]
 Taking particulary $f(t)=t^\beta$ we obtain:
  \[(\ln^\beta x)^{(n-k)}=x^{-n+k}\sum_{i=1}^{n-k}s(n-k,i)(\beta)_i\ln^{\beta-i}x.\]
Replacing these in (\ref{e3}) we have the following.
\begin{theorem} Let $\alpha\not=0$ be real, and $n,$  $i,\;(i\leq n)$ be nonnegative integers. Then
\begin{equation}\label{e4}s(n,i,\alpha)=\sum_{k=0}^{n-i}{n\choose k}(-\alpha)_ks(n-k,i).\end{equation}
\end{theorem}

\begin{note} Theorem is true even in the case $\alpha=0$ if we additionally define $(0)_0=1.$ \end{note}

\section{ Some combinatorial identities}

 Taking $i=1$ in (\ref{e4}) we obtain the following:
\begin{corollary}
Let $\alpha$ be a real number, and $n$ be a  positive  integer. Then
\[s(n,1,\alpha)=n!\sum_{k=0}^{n-1}(-1)^{n-k-1}\frac{{-\alpha\choose k}}{n-k}.\]
\end{corollary}
For $s(n,1,\alpha)$ we have the following recurrence relation:
\begin{equation}\label{e5} s(1,1,\alpha)=1,\;
s(n,1,\alpha)=(-\alpha-n+1)s(n-1,1,\alpha)+(-\alpha)_{n-1},\;(n\geq 2).\end{equation}

Particulary, we have $s(2,1,\alpha)=-2\alpha-1.$

We shall now prove that
polynomials $r(n,1,\alpha),\;(n=1,2,\ldots)$ defined by:
\[r(n,1,\alpha)=\sum_{k=0}^{n-1}(k+1)s(n,k+1)(-\alpha)^k,\]
 satisfy the above  recurrence relation. For $n=1$ it is obviously true.

 Using two terms recurrence relations for  Stirling numbers of the first kind  we obtain:
\[r(n,1,\alpha)=\sum_{k=0}^{n-1}(k+1)[s(n-1,k)-(n-1)s(n-1,k+1)](-\alpha)^k=\]\[=\sum_{k=0}^{n-1}(k+1)s(n-1,k)(-\alpha)^k-(n-1)\sum_{k=0}^{n-2}(k+1)s(n-1,k+1)(-\alpha)^k.\]

Since $s(n-1,0)=0,$ by
replacing $k+1$ instead of $k$ in the first sum on the right we have:
\[r(n,-\alpha)=\sum_{k=0}^{n-2}(k+2)s(n-1,k+1)(-\alpha)^{k+1}-(n-1)\sum_{k=0}^{n-2}(k+1)s(n-1,k+1)(-\alpha)^k=\]\[=
(-\alpha-n+1)r(n-1,\alpha)+\sum_{k=0}^{n-2}s(n-1,k+1)(-\alpha)^{k+1}.\]
Furthermore, a well known property of Stirling numbers implies:
\[ \sum_{k=0}^{n-2}s(n-1,k+1)(-\alpha)^{k+1}=\sum_{k=1}^{n-1}s(n-1,k)(-\alpha)^k=(-\alpha)_{n-1}.\]
We thus obtain that $r(n,1,\alpha)$ satisfies (\ref{e5}).
In this way we  have proved the following identity:
\begin{theorem} Let $\alpha$ be a real number, and $n\geq 1$ be an integer. Then:
\begin{equation}\label{e6}n!\sum_{k=0}^{n-1}(-1)^{k}\frac{{-\alpha\choose k}}{n-k}=\sum_{k=0}^{n-1}(k+1)\mathbf s(n,k+1)\alpha^k,\end{equation}
where $\mathbf s(n,k+1)$ are unsigned Stirling numbers of the firs kind.
\end{theorem}
\begin{note} Theorem is true for $\alpha=0$ with the convention that  $0^0=1.$\end{note}
Some particular values for $\alpha$  in (\ref{e6}) gives several interesting combinatorial identities.
For $\alpha=-1$ we obtain an identity expressing factorials in terms of Stirling numbers of the first kind.
\begin{corollary} The following formula is true:
$$(-1)^n(n-2)!=\sum_{k=0}^{n-1}(k+1)s(n,k+1),\;(n\geq 2).$$
\end{corollary}
For $\alpha=1$ we obtain a formula for the sum of reciprocals of natural numbers.
\begin{corollary} For each $n=1,2,\ldots$ we have:
$$n!\left(1+\frac 12+\cdots+\frac 1n\right)=\sum_{k=0}^{n-1}(k+1)\mathbf s(n,k+1).$$
\end{corollary}

 Consider now the case that $\alpha<0$ is an integer  and $n\geq 1-\alpha.$  In this case the function \[q(n,1,\alpha)=(-1)^{n+\alpha-1}(-\alpha)!(n+\alpha-1)!,\;(n\geq 1-\alpha)\] satisfies the recurrence relation (\ref{e5}).
In fact, in this case we have $(-\alpha)_{n-1}=0,$ since a factor in this product must be zero. Thus:
\[(-\alpha-n+1)q(n-1,1,\alpha)=-(n+\alpha-1)(-1)^{n-\alpha-2}(-\alpha)!(n+\alpha-2)!=\]\[=(-1)^{n-\alpha-1}(-\alpha)!(n+\alpha-1)!=q(n,1,\alpha).\]
 We have thus obtained the following:
\begin{corollary}
Let $n,\alpha$ be positive integers and $n\geq \alpha+1.$ Then:
\[n!\sum_{k=0}^{\alpha}\frac{(-1)^{\alpha-k}{\alpha\choose k}}{n-k}=\alpha!(n-\alpha-1)!.\]
In other wards, the following equation holds:
\[ (\alpha+1)\sum_{k=0}^{\alpha}\frac{(-1)^{\alpha-k}{\alpha\choose k}}{n-k}=\frac{1}{{n\choose \alpha+1}}.\]
\end{corollary}
In the  case that $\alpha$ is a negative integer and  $n\leq -\alpha$ the identity (\ref{e6}) is closely related to the harmonic numbers. Namely, denote $H_0=0,\;H_n=\sum_{k=1}^n\frac 1k,\;(n=1,2,\ldots),$ and define:
\[h(n,1,\alpha)=(H_{-\alpha}-H_{-\alpha-n})\frac{(-\alpha)!}{(-\alpha-n)!},\;(\alpha=-1,-2,\ldots;n=1,2,\ldots,-\alpha).\]
It will be shown that $h(n,1,\alpha)$ satisfies (\ref{e5}).
We use the induction on $n.$ For $n=1$ we have:
\[h(1,1,\alpha)=\frac{(-\alpha)!}{(-\alpha-1)!}(H_{-\alpha}-H_{-\alpha-1})=1.\]
Furthermore we have:  
\[(-\alpha-n+1)h(n-1,1,\alpha)+(-\alpha)_{n-1}=\frac{(-\alpha)!}{(-\alpha-n)!}(H_{-\alpha}-H_{-\alpha-n+1})+(-\alpha)_{n-1}=\]\[=
\frac{(-\alpha)!}{(-\alpha-n)!}(H_{-\alpha}-H_{-\alpha-n})-\frac{(-\alpha)!}{(-\alpha-n)!(-\alpha-n+1)}+(-\alpha)_{n-1}=\]\[=\frac{(-\alpha)!}{(-\alpha-n)!}(H_{-\alpha}-H_{-\alpha-n}),\]
since $\frac{(-\alpha)!}{(-\alpha-n)!(-\alpha-n+1)}=(-\alpha)_{n-1}.$
We thus obtain:
\[h(n,1,\alpha)=(-\alpha-n+1)b(n-1,1,\alpha)+(-\alpha)_{n-1},\] that is, $h(n,1,\alpha)$ satisfies (\ref{e5}). We have thus proved the following.
\begin{corollary} Let $\alpha$ be a positive integer and $1\leq n\leq \alpha$ be integers. Then:
\[H_\alpha-H_{\alpha-n}=\frac{(-1)^{n+1}}{{\alpha\choose n}}\sum_{k=0}^{n-1}\frac{(-1)^k{\alpha\choose k}}{n-k}.\]
Equivalently:
\[H_\alpha-H_{\alpha-n}=\frac{\sum_{k=0}^{n-1}(k+1)s(n,k+1)\alpha^k.}{\sum_{k=0}^{n}s(n,k)\alpha^k}.\]
\end{corollary}
In the case $\alpha=n$ we have the following expressions for harmonic numbers.
\begin{corollary} Let $n$ be a positive integer. Then:
\[H_n=(-1)^{n+1}\sum_{k=0}^{n-1}\frac{(-1)^k{n\choose k}}{n-k},\]
and
\[H_n=\frac{1}{n!}\sum_{k=0}^{n-1}(k+1)s(n,k+1)n^k.\]
\end{corollary}

\end{document}